\documentclass{article}
\usepackage{amssymb,amsmath,eucal}

\input{epsf}

\def\frV{\mathfrak V}

\def\frg{\mathfrak g}
\def\frh{\mathfrak h}

\def\frr{\mathfrak r}

\def\frv{\mathfrak v}

\def\cK{\mathcal K}
\def\cV{\mathcal V}

\def\Mor{\mathrm{Mor}}

\def\K{{\Bbb K}}

\def\Z{{\Bbb Z}}

\def\la{\langle}
\def\ra{\rangle}

\def\ov{\overline}
\def\phi{\varphi}
\def\epsilon{\varepsilon}
\def\kappa{\varkappa}
\def\le{\leqslant}
\def\ge{\geqslant}

\def\PBor{\mathrm{PsBor}}

\newcounter{sec}

\newcounter{punct}[sec]

\def\punct{\refstepcounter{punct}{\arabic{sec}.\arabic{punct}.  }}

\newtheorem{theorem}{Theorem}[sec]
\newtheorem{proposition}[theorem]{Proposition}

\def\G{\mathbb{G}}

\def\N{\mathbb{N}}

\def\COUNTERS{\addtocounter{sec}{1}
              \setcounter{punct}{0}
          \setcounter{equation}{0}
          \setcounter{theorem}{0}
          \setcounter{problem}{0}
          }
          
          \def\sm{\smallskip}
          \def\ov{\overline}
          \def \wt{\widetilde}
         \def \cO{\mathcal{O}}
\def\1{\mathbf 1}
          
\begin{document}

\begin{center}
{\bf\Large
Infinite symmetric group and bordisms of pseudomanifolds}

\bigskip

\sc \large  Alexander A. Gaifullin
\footnote{Supported by RFBR (project 14-01-92612), by grants of the President of 
 the Russian Federation (NSh-4833.2014.1 and MD-2969.2014.1), and by 
 Dmitri Zimin's ``Dynasty'' foundation.},
Yury A. Neretin%
\footnote{Supported by grants FWF, P22122, P25142}
\end{center}

{\small We consider a category whose morphisms are bordisms of $n$-dimensional
pseudomanifolds equipped with a certain additional structure (coloring). 
On the other hand,  we consider the  product $G$ of $(n+1)$ copies of infinite 
symmetric group. We show that unitary representations of $G$ produce functors
from the category of $(n-1)$-dimensional bordisms to the category of Hilbert spaces and
 bounded linear
operators.}

\section{Pseudomanifolds and pseudobordisms}

\COUNTERS

 First, we fix several definitions.
 
 \sm

 {\bf\punct Simplcial cell complexes.} Consider a disjoint union
$\coprod \Xi_j$ of a finite collection of simplices
$\Xi_j$. We consider a topological quotient space $\Sigma$ of 
$\coprod \Xi_j$  with respect to certain equivalence relation.
The quotient 
must satisfy the following properties

\sm

a) For any simplex $\Xi_i$, the tautological map $\xi_i:\Xi_i\to \Sigma$  is an embedding.
Therefore 
we can think of $\Xi_i$  as of a subset of $\Sigma$.

\sm

b) For any pair of simplices $\Xi_i$, $\Xi_j$, the intersection 
$\xi_i^{-1}\bigl(\xi_i(\Xi_i)\cap \xi_j(\Xi_j)\bigr)\subset \Xi_i$ is a union of faces of $\Xi_i$ and
the partially defined map
$$
\Xi_i\stackrel{\xi_i}{\longrightarrow} \Sigma\stackrel{\xi_j^{-1}}{\longrightarrow} \Xi_j
$$
is affine on each face.

\sm

We shall call such quotients  {\it simplicial cell complexes}.

\sm

{\sc Remark on terminology.} There are two similar (and more common) definitions of spaces composed from simplices
(see, e.g., \cite{Hat}). 
The first one is a more restrictive definition of ``a simplicial complex''. 
In this case, a non-empty intersection of two faces  is a (unique) face.
See examples of simplicial cell complexes, which are not simplicial complexes 
in Fig.2 and Fig.3.b.
A more wide class of simplicial spaces are $\Delta$-complexes, in this case 
glueing of a simplex with itself along faces is allowed (as for standard 1-vertex triangulations
of two-dimensional surfaces), see Fig 1.
\hfill $\lozenge$

\sm

{\bf\punct Pseudomanifolds.} A {\it pseudomanifold} of dimension $n$
is a simplicial cell complex
such that 

\sm

a) Each face is contained in an $n$-dimensional face. We call $n$-dimensional
faces  {\it chambers.}

\sm

b) Each $(n-1)$-dimensional face
is contained in precisely two chambers.

\sm

See, e.g., \cite{ST}, \cite{Gai1}.

\begin{figure}
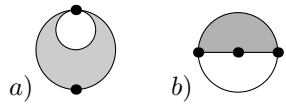

 $$
a) \epsfbox{draw1.3}\qquad  b)\epsfbox{draw1.4}
$$
\caption{To the definition of simplicial cell complexes. The triangle a) is forbidden, the pair of triangles b)
is allowed.}
\end{figure}

\sm

{\sc Remark.} Any  cycle of singular $\Z$-homologies in a topological space
can be realized as an image of a  pseudo-manifold
(this is more-or-less obvious).
Recall that  there are cycles in manifolds, which cannot be realized as  images of  manifolds.
\hfill $\lozenge$

\sm

{\sc Remark on terminology.} 
In literature, there exists  another variant of a definition
of a pseudomanifold. Seifert,  Threlfall, \cite{ST} impose two additional requiments:
a pseudomanifold must be a simplicial complex and must be 'strongly connected'.
The latter conditions means that the complement of the union of faces of codimension
2 must be connected.
\hfill $\lozenge$

\sm


\sm

{\bf \punct Normal pseudomanifolds and normalization.%
\label{ss:normal}}

{\it Links.} Let $\Sigma$ be a pseudomanifold,
let $\Gamma$ be its $k$-dimensional face. Consider all
$(k+1)$-dimensional faces $\Phi_j$ of $\Sigma$
containing $\Gamma$ and choose a point $\phi_j$
in the relative interior of each  face $\Phi_j$. For each face $\Psi_m\supset \Gamma$
we consider the convex hull of all points $\phi_j$ that are contained in $\Psi_m$. The link
of $\Gamma$ is the simplicial cell complex whose faces are such convex hulls.

\sm

{\it Normal pseudomanifolds.}
 A pseudomanifold is {\it normal} if the  link of any face of codimension
 $\ge 2$
is connected.

\sm

{\sc Example.} Consider a triangulated compact two-dimensional surface $\Sigma$.
Let $a$, $b$ be two vertices that are not connected by an edge. Glueing together
$a$ and $b$ we get a pseudomanifold which is not normal, see Fig.2. 
\hfill $\lozenge$

\sm

{\it Normalization.}
For any pseudomanifold $\Sigma$ there is a unique {\it normalization}
(\cite{GM}),
i.e. a normal pseudomanifold $\wt\Sigma$ and a map $\pi:\wt\Sigma \to \Sigma$ such that

\sm

--- restriction of $\pi$ to any face of $\wt\Sigma$ is an affine bijective map
of faces.

\sm

--- the map $\pi$ send different $n$-dimensional and $(n-1)$-dimensional 
faces to different faces. 

\sm

\begin{figure}
$$
 \epsfbox{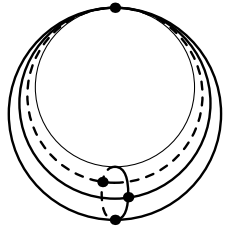}
$$

\caption{ A non-normal two-dimensional pseudomanifold.}
\end{figure}

{\it A construction of the normalization.}
To obtain a normalization of $\Sigma$ we cut a pseudomanifold $\Sigma$ 
into a disjoint collection of chambers
$\Xi_i$. As above, denote by $\xi_i:\Xi_i\to \Sigma_j$ the embedding of $\Xi_i$ to $\Sigma$.
Let $x\in\Xi_i$, $y\in \Xi_j$. We say that $x\sim y$ if $\xi_i(x)=\xi_j(y)$ and 
this point is contained in a common $(n-1)$-dimensional face of the chambers $\xi_i(\Xi_i)$
and $\xi_j(\Xi_j)$. We extend $\sim$ to an equivalence relation by the transitivity.
The quotient of $\coprod \Xi_i$ is the normalization of $\Sigma$.

The following way of normalization  is more visual. Let $\Sigma$ be non-normal.
 Let $\Xi$ be a face
of codimension $2$ with
 link  consisting of $m$ connected components. 
Consider a small closed neighborhood $\cO$ of $\Xi$ in $\Sigma$. Then $\cO\setminus \Xi$ is disconnected
and consists of $m$ components, say $\cO_1$,\dots, $\cO_m$.
Let $\ov\cO_j$ be the closure of $\cO_j$ in $\Sigma$, $\ov\cO_j=\cO_j\cup \Xi_j$.
We replace $\cO$ by the disjoint union of $\ov\cO_j$ and get a new pseudomanifold
$\Sigma'$ (in Fig.2, we duplicate the upper vertex). Then we repeat the same operation to another 
stratum with disconnected link. These operation enlarges
number of strata of codimension $\ge 2$, the strata of dimension $n$ and $(n-1)$ remain
the same (and the incidence of these strata is preserved). Therefore the process is finite and we get a normal pseudomanifold.
  \hfill $\lozenge$

\smallskip

{\bf\punct Colored pseudomanifolds.}
Consider an $n$-dimensional {\it normal} pseudomanifold
$\Sigma$. A coloring of $\Sigma$  is the following structure

\sm

a) To any chamber we assign a sign $(+)$ or $(-)$. Chambers adjacent to plus-chambers are 
minus-chamber and vise versa.

\sm

b) Choose $n+1$ colors (say, red, blue, green,  orange, etc.).
 Each vertex of the complex is colored in such a way that
the colors of vertices of each chamber are pairwise different.

\sm

c) All $(n-1)$-dimensional faces are colored, 
in such a way that colors of faces of a chamber
are pairwise different, and a color of a face coincides with a color of the opposite
vertex of any chamber containing this face.

\begin{figure}
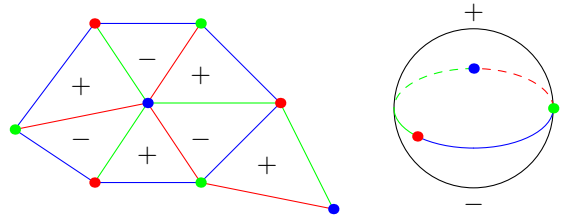

$$\epsfbox{draw1.1}\qquad \epsfbox{draw1.2}$$

\caption{Reference to the definition of colored pseudomanifolds:
\newline
a) a colored two-dimensional pseudomanifold;
\newline
b) a double chamber.}

\end{figure}

\sm

 We say that a {\it double-chamber} is a colored $n$-dimensional
 pseudomanifold obtained from two identical
 copies $\Delta_1$, $\Delta_2$ of an $n$-dimensional simplex by identification of the corresponding
 $x\in \Delta_1$, $x\in\Delta_2$ of the boundaries of $\Delta_1$, $\Delta_2$.
 
 \sm
 
{\sc Remark.} Colored pseudomanifolds were introduced by Pezzana
and Ferri in 1975-1976, see  \cite{Pez}, \cite{Fer}, \cite{FGG}, \cite{Gai}.
\hfill $\lozenge$

\sm

{\bf\punct Colored pseudobordisms.%
\label{ss:pseudobordisms}} Fix $n\ge 1$.
We define a category $\PBor$ of pseudobordisms.
Its objects are nonnegative integers. 
A morphism $\beta\to\alpha$ is
the following collection of data

\sm

1) A colored $n$-dimensional pseudomanifold (generally, disconnected).

\sm

2) An injective map of the set $\{1,2,\dots,\alpha\}$ to the set of plus-chambers
and an injective map of the set $\{1,2,\dots,\beta\}$ to the set of minus-chambers
In other words, we assign labels $1$, \dots, $\alpha$ to some plus-chambers.
and labels $1$, \dots, $\beta$ to some minus-chambers.

\sm

We require that each double-chamber has at least one label.

\smallskip

{\it Composition.}
Let $\Sigma\in \Mor(\beta,\alpha)$, $\Lambda\in\Mor(\gamma,\beta)$.
We define their composition $\Sigma\diamond \Lambda$ as follows.
Remove interiors of labeled minus-chambers of $\Sigma$
and interiors of
 labeled plus-chambers of $\Lambda$. Next, for each 
 $s\le \beta$, we glue boundaries of the minus-chamber of $\Sigma$
 with label $s$ with the boundary of the plus-chamber of $\Lambda$ with label
 $s$ according the simplicial structure of boundaries  and colorings of
 $(n-1)$-simplices. Next, we normalize the resulting pseudomanifold.
 
Finally we remove label-less double chambers (such components can arise as a result of
gluing of two label-keeping double chambers).
 
 \sm
 
 {\it Involution.}
 For a morphism $\Sigma\in  \Mor(\beta,\alpha)$ 
 we define the morphism $\Sigma^*\in \Mor(\alpha,\beta)$
 by changing of signs on chambers. Thus we get an {\it involution}
 in the category $\PBor$. For any
 $T\in \Mor(\beta,\alpha)$, $S\in\Mor(\gamma,\beta)$
 we have 
 $$
 (S\diamond T)^*=T^*\diamond S^*
 $$


\sm

{\bf \punct Further structure of the paper.} Below we construct a family of functors
from the category of pseudobordisms to the category of Hilbert spaces and bounded operators.
This means that for each $\alpha$ we construct a Hilbert space $H(\alpha)$
and for each morphism $\Sigma\in\Mor(\beta,\alpha)$ we construct an operator
$\rho(\Sigma):\,H(\beta)\to H(\alpha)$ such that
for any $\alpha$, $\beta$, $\gamma$ and 
any $\Sigma\in \Mor(\beta,\alpha)$, $\Xi\in\Mor(\gamma,\beta)$,
$$
\rho(\Sigma\diamond \Xi)=\rho(\Sigma)\, \rho(\Xi) 
.$$
Also, representations obtained below satisfy properties
$$
\rho(\Sigma)^*=\rho(\Sigma^*),\qquad \|\rho(\Sigma)\|\le 1
.
$$ 
In fact, such functors  arise in a natural way
from the representation theory of infinite symmetric groups. 
In Section 2, we introduce a category of double cosets
on the product of $(n+1)$ copies of an infinite symmetric group.
In Section 3, we show that
the category of double cosets  is equivalent to the category of pseudo-bordisms
(this is the main statement of this note).
In Section 4, we construct a family of representations of this category
(statements of this section are more-or-less automatic).

For $n=1$ the construction  reduces to Olshanski \cite{Olsh1},
 for $n=2$ it coincides with \cite{Ner-tri}.  


\section{Multiplication of double cosets}

\COUNTERS

{\bf\punct Symmetric groups. Notation. } Let $S(L)$ be the group of permutations
of a set with $L$ elements.
Denote by $K=S(\infty)$ the group of finitely supported permutations of 
$\N$. By $\ov K=\ov S(\infty)$ we denote the group of all permutations of $\N$.
Denote by $K(\alpha)\subset K$, $\ov K(\alpha)\subset \ov K$ the stabilizers of points
$1$, \dots, $\alpha$. 
We equip $\ov S(\infty)$ with a natural topology assuming that
the subgroups $K(\alpha)$ are open.

Sometimes we will represent elements of symmetric groups as $0-1$-matrices.

\smallskip


{\bf\punct Multiplication of double cosets.} Denote the product
of $(n+1)$ copies of  $S(\infty)$ 
by $G$. Denote by $K\simeq S(\infty)$ the diagonal subgroup in $G$,
its elements have the form $(g,g,\dots,g)$.


Consider double cosets 
$K(\alpha)\setminus G/K(\beta)$, i.e., elements of $G$ defined up to the equivalence
$$
g\sim k_1 g k_2,\qquad k_1\in K(\alpha),\, k_2\in K(\beta)
$$
We wish to define product of double cosets
$$
K(\alpha)\setminus G/K(\beta)\,\,\times\,\, K(\beta)\setminus G/K(\gamma)\,\,\to\,\,
K(\alpha)\setminus G/K(\gamma).
$$

For this purpose, define elements $\theta_\sigma[j]\in K(\sigma)$ by
$$
\theta_\sigma[j]:=
\begin{pmatrix}
 \1_\sigma&0&0&0
 \\
 0&0&\1_j&0
 \\
 0&\1_j&0&0
 \\
 0&0&0&\1_\infty
\end{pmatrix}
,$$
where $\1_j$ denotes the unit matrix of order $j$.

\begin{proposition}
 Let 
 $$\frg\in K(\alpha)\setminus G/ K(\beta), \quad\frh\in K(\beta)\setminus G/K(\gamma)$$
 be double cosets. Let $g$, $h\in G$ be their representatives.
 Then the sequence
 \begin{equation}
\frr_j:= K(\alpha)\cdot g \theta_\beta[j] h\cdot K(\gamma)\,\,\,\in K(\alpha)\setminus G/K(\gamma)
 \label{eq:theta-1}
 \end{equation}
 is eventually constant.
 The limit value of $\frr_j$ does not depend on a choice of representatives $g\in\frg$ and $h\in\frh$. 
 Moreover, if $g$, $h\in S(L)^{n+1}\subset S(\infty)^{n+1}$,
 then it is sufficient to consider $j=L-\beta$.
\end{proposition}

We define the product 
$$\frg \circ \frh\in K(\alpha)\setminus G/K(\gamma)$$
 of double cosets as
the limit value of the sequence (\ref{eq:theta-1}). 

\smallskip

{\sc Formula for product.}
 Represent $g$ as a collection of block matrices $\bigl(g^{(1)},\dots, g^{(n+1)}\bigr)$ of size
$$\bigl(\alpha+(L-\alpha)+(L-\beta)+\infty\bigr)\times
\bigl(\beta+(L-\beta)+(L-\beta)+\infty\bigr),$$
represent $h$
as a  collection of block matrices $\bigl(h^{(1)},\dots, h^{(n+1)}\bigr)$ of size 
$$\bigl(\beta+(L-\beta)+(L-\beta)+\infty\bigr)\times\bigl(\gamma+(L-\gamma)+(L-\beta)+\infty\bigr)$$
\begin{equation}
g^{(k)}=
\begin{pmatrix}
a^{(k)}&b^{(k)}&0&0\\c^{(k)}&d^{(k)}&0&0\\
0&0&\1_{L-\beta}&0
\\0&0&0&\1_\infty
\end{pmatrix}
,\qquad
h^{(k)}=
\begin{pmatrix}
p^{(k)}&q^{(k)}&0&0\\r^{(k)}&t^{(k)}&0&0
\\0&0&\1_{L-\beta}&0\\0&0&0&\1_\infty
\end{pmatrix}
.
\label{eq:gh}
\end{equation}
Then we write a representative of the double coset $\frg\circ\frh$ as
$$
(g\circ h)^{(k)}:= g\cdot \theta_{\beta}[L-\beta]\cdot h
= 
\begin{pmatrix} a^{(k)}p^{(k)}&a^{(k)}q^{(k)}&b^{(k)}&0\\
c^{(k)}p^{(k)}&c^{(k)}q^{(k)}&d^{(k)}&0\\r^{(k)}&t^{(k)}&0&0\\
0&0&0&\1_\infty
 \end{pmatrix}
.
$$

{\sc Proof.} First, we show that the result does not depend on a choice of $j$
 Denote 
$$\mu=L-\beta, \quad \nu=L-\alpha,\quad \kappa=L-\gamma
.$$
 Preserving the previous notation for
$g^{(k)}$, $h^{(k)}$, we write
$$
(g\cdot \theta_\beta[\mu+j]\cdot  h)^{(k)}=
\begin{pmatrix}
a^{(k)}p^{(k)}&a^{(k)}q^{(k)}&0&b^{(k)}&0&0\\
c^{(k)}p^{(k)}&c^{(k)}q^{(k)}&0&d^{(k)}&0&0\\
0&0&0&0&\1_j&0\\
r^{(k)}&t^{(k)}&0&0&0&0\\
0&0&\1_j&0&0&0\\
0&0&0&0&0&\1_\infty
\end{pmatrix}
.$$
This coincides with 
\begin{multline*}
\begin{pmatrix}
\1_\alpha&0&0&0&0&0\\
0&\1_\nu&0&0&0&0\\
0&0&0&\1_j&   0&0\\
0&  0&\1_\mu &0&0&0\\
0&0&0&0&\1_j&0\\
0&0&0&0&0&\1_\infty
\end{pmatrix}
\begin{pmatrix} a^{(k)}p^{(k)}&a^{(k)}q^{(k)}&b^{(k)}&0&0&0\\
c^{(k)}p^{(k)}&c^{(k)}q^{(k)}&d^{(k)}&0&0&0
\\r^{(k)}&t^{(k)}&0&0&0&0\\
0&0&0&\1_j&0&0\\
0&0&0&0&\1_j&0\\
0&0&0&0&0&\1_\infty
 \end{pmatrix}
 \times \\ \times
\begin{pmatrix}
\1_\gamma&0&0&0&0&0\\
0&\1_\kappa&0&0&0&0\\
0&0&0&\1_\mu&   0&0\\
0&0&0&0&\1_j&0\\
0&  0&\1_j &0&0&0\\
0&0&0&0&0&\1_\infty
\end{pmatrix}
.
\end{multline*}

Next, we show that  (\ref{eq:theta-1}) does not depend on the
choice of representatives of double cosets. To be definite, replace
a collection $\{g^{(k)}\}$ in (\ref{eq:gh}) by 
$$
\begin{pmatrix}
\1_\alpha &0&0\\0&u&0\\0&0&\1_\infty
\end{pmatrix}
\begin{pmatrix}
a^{(k)}&b^{(k)}&0\\c^{(k)}&d^{(k)}&0\\0&0&\1_\infty
\end{pmatrix}
\begin{pmatrix}
\1_\alpha &0&0\\0&v&0\\0&0&\1_\infty
\end{pmatrix}
=
\begin{pmatrix}
a^{(k)}&b^{(k)}v&0\\uc^{(k)}&ud^{(k)}v&0\\0&0&\1_\infty
\end{pmatrix}
.
$$
Then $(g\circ h)^{(k)}$ is
\begin{multline*}
\begin{pmatrix} a^{(k)}p^{(k)}&a^{(k)}q^{(k)}&b^{(k)}v&0\\
uc^{(k)}p^{(k)}&uc^{(k)}q^{(k)}&ud^{(k)}v&0\\r^{(k)}&t^{(k)}&0&0\\
0&0&0&\1_\infty
 \end{pmatrix}
 =\\=
  \begin{pmatrix} \1_\alpha&0&0&0\\
0&u&0&0\\0&0&\1_\mu &0\\
0&0&0&\1_\infty
\end{pmatrix}
 \begin{pmatrix} a^{(k)}p^{(k)}&a^{(k)}q^{(k)}&b^{(k)}&0\\
c^{(k)}p^{(k)}&c^{(k)}q^{(k)}&d^{(k)}&0\\r^{(k)}&t^{(k)}&0&0\\
0&0&0&\1_\infty
 \end{pmatrix}
   \begin{pmatrix} \1_\gamma&0&0&0\\
0&\1_\kappa &0&0\\0&0&v &0\\
0&0&0&\1_\infty
\end{pmatrix}
\end{multline*}
This completes the proof.
\hfill $\square$

\smallskip

\begin{proposition}
 The $\circ$-product is associative.
\end{proposition}

{\sc Proof.} Let $g$, $h\in G$ be  as above,
and  let $w=(w^{(1)},\dots, w^{(n+1)})\in G$ be given by
$$
w^{(k)}=
\begin{pmatrix}
x^{(k)}&z^{(k)}&0\\y^{(k)}&u^{(k)}&0\\0&0&\1_\infty
\end{pmatrix}
.
$$
Evaluating $(g\circ h)\circ w$ and $g\circ (h\circ w)$
we get the matrices
\begin{equation}
\begin{pmatrix}
a^{(k)}p^{(k)}x^{(k)}& a^{(k)}p^{(k)}y^{(k)} & a^{(k)}q^{(k)} & b^{(k)}& 0& 0\\
c^{(k)}p^{(k)}x^{(k)}& c^{(k)}p^{(k)}y^{(k)} & c^{(k)}q^{(k)} & d^{(k)}& 0& 0\\
0&0&0&0&\1&0\\
r^{(k)}x^{(k)} & r^{(k)}y^{(k)} & t^{(k)} &0&0& 0\\
z^{(k)}&u^{(k)}&0&0&0& 0\\
0&0&0&0& 0&\1_\infty
\end{pmatrix}
\label{eq:matrix-1}
,
\end{equation}
\begin{equation}
\begin{pmatrix}
a^{(k)}p^{(k)}x^{(k)}& a^{(k)}p^{(k)}y^{(k)} & a^{(k)}q^{(k)}& 0 & b^{(k)}& 0\\
c^{(k)}p^{(k)}x^{(k)}& c^{(k)}p^{(k)}y^{(k)} & c^{(k)}q^{(k)}& 0 & d^{(k)}& 0\\
r^{(k)}x^{(k)} & r^{(k)}y^{(k)} & t^{(k)}& 0& 0 &0\\
z^{(k)}&u^{(k)}&0&0&0& 0\\
0&0&0& \1&0& 0\\
0&0&0&0& 0&\1_\infty
\end{pmatrix}
\label{eq:matrix-2}
.
\end{equation}
Both matrices are elements of the double coset containing
\begin{equation}
\begin{pmatrix}
a^{(k)}p^{(k)}x^{(k)}& a^{(k)}p^{(k)}y^{(k)} & a^{(k)}q^{(k)} & b^{(k)}& 0& 0\\
c^{(k)}p^{(k)}x^{(k)}& c^{(k)}p^{(k)}y^{(k)} & c^{(k)}q^{(k)} & d^{(k)}& 0& 0\\
r^{(k)}x^{(k)} & r^{(k)}y^{(k)} & t^{(k)} &0&0& 0\\
z^{(k)}&u^{(k)}&0&0&0& 0\\
0&0&0&0& \1&0\\
0&0&0&0& 0&\1_\infty
\end{pmatrix}
\label{eq:matrix-3}
,\end{equation}
matrix (\ref{eq:matrix-1}) is obtained from (\ref{eq:matrix-3})
by a permutation of rows, and matrix (\ref{eq:matrix-2})
is obtained from (\ref{eq:matrix-3}) by a permutation of columns.
\hfill $\square$

\smallskip

Thus we get a category $\cK$, whose objects are nonnegative integers,
and
$$
\Mor_\cK(\beta,\alpha):=K(\alpha)\setminus G/K(\beta)
.
$$
A product is the product of double cosets.

\smallskip

{\bf\punct Involution.} The map
$g\mapsto g^{-1}$ induces the map $\frg\mapsto \frg^*$ of double cosets
$$
K(\alpha)\setminus G/K(\beta) \to K(\beta)\setminus G/K(\alpha)
.$$
Evidently, $(\frg\circ \frh)^*=\frh^*\circ \frg^*$.


\section{Correspondence}

\COUNTERS

{\bf\punct A correspondence between pseudomanifolds and symmetric groups.%
\label{ss:correspondence}} Denote by $S(L)$ the symmetric group of order $L$.
 Denote by 
$$
S(L)^{n+1}:=S(L)\times\dots\times S(L)
$$
the direct product of $n+1$ copies of $S(L)$, we assign $n+1$ colors,
say, red,  blue,  orange, etc.,
to copies of $S(L)$.

Consider a colored pseudomanifold $\Sigma$ with $2L$ chambers. We say that
a labeling of $\Sigma$ is a bijection of the set $\{1,2,\dots,L\}$ with the set of
plus-chambers of $\Sigma$ and a bijection of  $\{1,2,\dots,L\}$ with the set of
minus-chambers of $\Sigma$.

\begin{theorem}
 There is a canonical one-to-one correspondence between the group $S(L)^{n+1}$ and 
 the set of all labeled colored normal
 $n$-dimensional pseudo-manifolds with $2L$ chambers.
\end{theorem}

{\sc Remark.} This correspondence for $n=2$ was proposed 
in \cite{Ner-tri}.
Earlier there was a construction of Pezzana--Ferri
(1975-1976), \cite{Pez}, \cite{Fer}, \cite{FGG}.
They considered  bipartie $(n+1)$-valent graphs whose edges are colored 
in $(n+1)$ colors, edges adjacent to a given vertex have pairwise different colors.
Such graphs correspond to colored pseudomanifolds. In \cite{Gai}--\cite{Gai3}
there was considered an action of free product $\Z_2*\dots*\Z_2$ of $n$ copies of $\Z_2$
on the set of chambers of a colored pseudomanifold. A construction relative
to the present construction was considered in \cite{BG}.
\hfill $\lozenge$

\sm

{\sc Construction of the correspondence.}
Indeed, consider a labeled colored normal pseudomanifold $\Sigma$ with $2L$ chambers.
Fix a color (say, blue). Consider all blue $(n-1)$-dimensional faces $A_1$, $A_2$, \dots.
Each blue face $A_j$ is contained in the plus-chamber with some label
$p(j)$ and in the minus-chamber with some label $q(j)$. We take an element
of the symmetric group $S(L)$ setting $p(j)\mapsto q(j)$ for all blue faces $A_j$.
We repeat the same construction
for all colors and obtain a tuple $(g^{(1)},\dots,g^{(n+1)})\in S(L)^{n+1}$.

Conversely, consider an element of the group $S(L)^{n+1}$. Consider $L$ labeled copies
of a colored chamber
(plus-chambers)
and another collection of $L$ labeled copies of the same chamber with another orientation
(minus-chambers). Let the blue permutation send $\alpha\mapsto\beta$. Then we glue the
the plus-chamber with label $\alpha$ with the minus-chamber  with label $\beta$ along the blue face
 (preserving colorings 
of vertices). The same is done for all colors.
 The  obtained  pseudomanifold $\Sigma$ is normal because the normalization
procedure from Subsection \ref{ss:normal} applied to $\Sigma$ produces $\Sigma$ itself.
\hfill $\square$

\smallskip

{\bf\punct The multiplication in symmetric group and pseudomanifolds.} 
Describe the multiplication in $S(L)^{n+1}$ in a geometric language.
Consider two labeled colored pseudomanifolds $\Sigma$, $\Xi$.
Remove interiors of minus-chambers of $\Sigma$ remembering a minus-label on each face of a removed
chamber, denote the topological space obtained in this way
by $\Sigma_-$. All $(n-1)$-faces of $\Sigma_-$ are colored and labeled.
In the same way, we remove plus-chambers from $\Xi$ and get a complex $\Xi_+$.
Next, we glue the corresponding faces of $\Sigma_-$ and $\Xi_+$ 
(with coinciding colors and labels according coloring
of vertices). In this way, we get a pseudomanifold and consider its normalization.

\smallskip

{\bf\punct Infinite case.} 
We say that an {\it infinite pseudo-manifold} is a disjoint union of a countable collection of
compact pseudomanifolds such that all but a finite number of
its components are double-chambers.

We define a colored infinite pseudo-manifold as above. A labeled pseudomanifold
is a colored  pseudomanifold with a numbering of plus-chambers by natural numbers and
a numbering of minus-chambers by natural
numbers such that all but a finite number of double-chambers have the same labels on both chambers.

\begin{theorem}
There is a canonical one-to-one correspondence between the group $S(\infty)^{n+1}$
and the set of all labeled colored normal infinite pseudomanifolds.
\end{theorem}

The correspondence is given by the same construction obtained  as above.

\sm

{\bf \punct Equivalence of categories.} 

\begin{theorem}
The category $\cK$ of double cosets and the category $\PBor$
of pseudo-bordisms are equivalent.
The equivalence is given by the following construction.
\end{theorem}

{\sc Correspondence $\Mor_\cK(\beta,\alpha)\longleftrightarrow \Mor_\PBor(\beta,\alpha)$.}
Let $\frg\in K(\alpha)\setminus G/ K(\beta)$ be a double coset. Let
$g\in \frg$ be its representative. Consider the corresponding labeled colored
pseudomanifold. A left multiplication $g\mapsto ug$ by an element
$u\in K(\alpha)$ is equivalent to a permutation $u$ of labels 
$\alpha+1$, $\alpha+2$, \dots on plus-chambers. A right multiplication $g\mapsto gv$ by an element $v\in K(\beta)$ is equivalent
 to a permutation of labels $\beta+1$, $\beta+2$, \dots on minus-chambers.
 
 Thus passing to double cosets is equivalent to forgetting  labels $>\alpha$ on
 plus-chambers and labels $>\beta$ on minus-chambers.
  Notice that all but a finite  number of double-chambers are label-less.
 Such label-less double chambers can be forgotten.
 Thus we get a pseudobordism.
 
\smallskip
 
{\sc Correspondence of products.} Let $g$, $h$ be representatives of double cosets.
Let $\Sigma$, $\Xi$ be the corresponding infinite labeled colored pseudomanifolds.
Let $\Sigma'$ correspond to $g \theta_\beta[j]$ , where $j$ is large.
We multiply $g \theta_\beta[j]$ by $h$ according to the rule
in Subsection \ref{ss:correspondence}.

Notice that minus-chambers of $\Sigma'$ with labels $>\beta$ are glued with 
double-chambers. Plus-chambers of $\Xi$ with labels $>\beta$
are also  glued with double chambers.
Both operations yield a changing of labels on chambers.
 This means
that in fact we glue together  only chambers with labels $\le\beta$,
in remaining cases we  change labels on chambers only. 
Afterwards we forget all labels which are grater than $\beta$  and get the operation
described in Subsection \ref{ss:pseudobordisms}.
\hfill $\square$


 \section{Representations}
 
 \COUNTERS

 Here we construct a family of representations
 of the group $G$. This produces representations of the category of double
 cosets and therefore representations of the category of pseudobordisms.
 The construction is a special case of \cite{Ner-tri} (where  the case $n=2$ was considered),
 more ways of constructions
 of representations of the group $G$, see in \cite{Ner-comb}, \cite{Ner-tri}.

\sm
 
 {\bf\punct The group $\G$.} We define an 'intermediate' group 
 $\G$,
 $$
 S(\infty)^{n+1}\subset \G\subset \ov S(\infty)^{n+1}
 ,$$
consisting of tuples $(g_1,\dots,g_{n+1})\in \ov S(\infty)^{n+1}$
such that $g_i g_j^{-1}\in S(\infty)$ for all $i$, $j$.
Denote by $\K\simeq \ov S(\infty)$ the diagonal subgroup consisting of 
tuples $(g,\dots,g)$. Define the subgroup $\K(\alpha)$ to be the group
of all $(h,\dots,h)$, where $h$ fixes $1$, \dots, $\alpha$.
Define the topology on $\G$ assuming that subgroups $\K(\alpha)$ are open.

 Obviously, there is the  identification
of double cosets
$$
K(\alpha)\setminus G/K(\beta)\simeq \K(\alpha)\setminus \G/\K(\beta)
.
$$

\smallskip
 
{\bf\punct A family of representation of $\G$.}  Consider $(n+1)$
Hilbert spaces%
\footnote{We admit arbitrary, finite-dimensional or infinite-dimensional,
separable Hilbert spaces.} $V_{red}$, $V_{orange}$, $V_{blue}$,  \dots.
Consider their tensor product
$$
\cV=V_{red}\otimes V_{blue} \otimes V_{green}\otimes \dots
$$
Fix a unit vector $\xi\in \cV$. Consider a countable tensor 
product of Hilbert spaces
\begin{multline}
\frV=
(\cV,\xi)\otimes (\cV,\xi)\otimes (\cV,\xi)\otimes \dots
=\\=
(V_{red}\otimes V_{blue} \otimes \dots\,\,,\xi)\otimes
(V_{red}\otimes V_{blue} \otimes\dots\,\,,\xi)\otimes
\dots
\label{eq:frV}
\end{multline}
 (for a definition of tensor products, see \cite{vN}). Denote 
 $$
 \frv=\xi\otimes \xi\otimes\dots \in \frV
 .
 $$
We define a representation $\nu$  of $\G$ 
in $\frV$ in the following way.
 The 'red' copy of $S(\infty)$ acts by
permutations of factors $V_{red}$. The 'blue' copy $S_\infty$ acts
by permutation of factors $V_{blue}$, etc. Thus we get an action of
the group $S(\infty)^{n+1}$.
The diagonal $\K=\ov S(\infty)$ acts by permutations of 
factors $\cV$. 

\sm

{\sc Remark.} For type I groups  $H_1$, $H_2$ irreducible unitary representations of 
$H_1\times H_2$ are tensor products of representations of $H_1$ and $H_2$  (see, e.g., \cite{Dix}
13.1.8). However,
$S(\infty)$ is not a type I group. {\it  Representations of $S(\infty)^{n+1}$ 
constructed above 
are not tensor products of representations of $S(\infty)$.}
\hfill $\lozenge$

\smallskip


{\bf \punct Representations of the category $\cK$.} Consider a unitary representation
$\rho$ of the group $\G$ in a Hilbert space $H$.
For $\alpha=0$, $1$, $2$, \dots consider the subspace $H_\alpha$
of $\K(\alpha)$-fixed vectors in $H$. Denote by $P_\alpha$ the operator 
of orthogonal projection to $H_\alpha$.
Let $\frg\in
 \K(\alpha)\setminus \G/\K(\beta)$ be a double coset,
 and let
$g\in \G$ be its representative. We define
an operator
$$
\ov\rho(\frg): H_\beta\to H_\alpha
$$
by 
$$
\ov\rho(\frg)=P_\alpha \rho(g)\Bigr|_{H_\beta}
$$

\begin{theorem}
\label{th:multiplicativity}
The operator $\ov\rho(\frg)$ does not depend on the choice
of a representative $g\in \frg$. For any $\alpha$, $\beta$, $\gamma$,
$$
\frg\in \K(\alpha)\setminus \G/\K(\beta),\quad \frh\in \K(\beta)\setminus \G/\K(\gamma)
$$
we have
$$
\ov \rho(\frg)\ov\rho(\frh)=\ov\rho(\frg\circ \frh)
$$
\end{theorem}

See a proof for  $n=2$ in \cite{Ner-tri},
the general case  is  completely similar (also this is a special case 
of \cite{Ner-cat}, Theorem VIII.5.1)

\begin{theorem}
Let $\pi$ be a representation  of the category $\cK$ in Hilbert spaces compatible with the involution
and satisfying $\|\pi(\frg)\|\le 1$ for all $\frg$.
Then $\pi$ is equivalent
to some representation
$\ov\rho$, where $\rho$ is a unitary representation of $\G$.
\end{theorem}

This is a special case of \cite{Ner-cat}, Theorem VIII.1.10.

\smallskip

{\bf\punct Spherical functions.}
In the  above example  we have
$$
\frV_\alpha=\underbrace{(\cV,\xi)\otimes \dots (\cV,\xi)}_{
\text{$\alpha$ times}}\otimes \xi\otimes \xi\dots\simeq
\cV^{\otimes \alpha}
,$$
in particular
$$
\frV_0=\frv.
$$

We wish to write an explicit formula for the spherical function
$$
\Phi(g)=\la \nu(g) \frv,\frv\ra.
$$

\begin{figure}
 $$
 \epsfbox{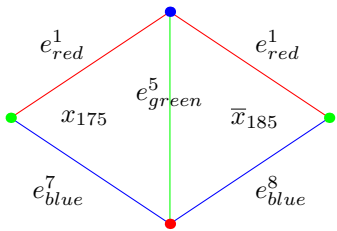}
 $$
 \caption{Arrangement of basis elements on a pseudomanifold}
\end{figure}

Choose an orthonormal basis in each space $V_{red}$, $V_{blue}$, $V_{green}$, etc.
$$
e_{red}^{i}\in V_{red}, \quad e_{blue}^j\in V_{blue},\quad e^k_{green}\in V_{green},\dots
$$
This determines the basis 
\begin{equation*}
e_{red}^{i}\otimes e_{blue}^j\otimes e^k_{green}\otimes \dots
\label{eq:vector}
\end{equation*}
in $\cV$. Expand $\xi$ in this basis,
\begin{equation}
\xi=\sum x_{ijk\dots} e_{red}^{i}\otimes e_{blue}^j\otimes e^k_{green}\otimes \dots
\label{eq:x}
\end{equation}

Consider the double coset $\frg$ containing $g$ and the corresponding
 colored pseudomanifold $\Sigma$.
Assign to each $(n-1)$-face an element of the basis of the corresponding color (in arbitrary way). 
Fix such arrangement. Consider a chamber $\Delta$, on its faces we have certain basis vectors
$e_{red}^{i}$, $e_{blue}^j$, $e^k_{green}$, ... .  Then we assign the number
$x(\Delta):=x_{ijk\dots}$ (see the last formula) to $\Delta$.

 \begin{proposition}
 $$
 \Phi(g)=
 \sum\limits_{\begin{matrix}\text{\scriptsize arangements}\\
\text{\scriptsize  of basis elements} \end{matrix}}
 \prod\limits_{\text{plus-chambers $\Delta$}} x(\Delta) \cdot
 \prod\limits_{\text{minus-chambers $\Gamma$}} \ov{x(\Gamma)}
 $$
 \end{proposition}
 
 Proof coincides with proof of Proposition 4.2 in \cite{Ner-tri}.

\tt
\noindent
A.Gaifullin: Steklov Mathematical Institute, Moscow, Russia;\\
MechMath Dept., Moscow State University\\
Kharkevich Institute for Information Transmission Problems, Moscow, Russia\\
agaif@mi.ras.ru

\noindent
Yu.Neretin: Math. Dept., University of Vienna; \\
Institute for Theoretical and Experimental Physics; \\
MechMath Dept., Moscow State University
\\
neretin(at)mccme.ru

\end{document}